\documentclass[12pt]{amsart}

\title{Free Diffusions and Property AO}
\author{Jason Asher}

\usepackage{amsmath}
\usepackage{amssymb}
\usepackage{fullpage}

\newtheorem{thm}{Theorem}[section]
\newtheorem{lem}[thm]{Lemma}

\def\cxmnsa{\mathbb{C}\langle X_1,\ldots,X_m,X_1^*,\ldots,X_m^*\rangle}
\def\cxm{\mathbb{C}\langle X_1,\ldots,X_m\rangle}

\begin{document}
\maketitle

\section{Introduction}

Guionnet and Shlyakhtenko extended Langevin-type free stochastic dynamics to the case of operators interacting by a locally convex potential in \cite{gushdiffusions}, and among other things used these results to give technical properties of certain operator algebras. They were specifically interested in algebras generated by the stationary laws of free stochastic differential equations (SDE) of the form
$$ dX_t = dS_t - \frac{1}{2}DV(X_t)dt $$
for a suitably locally convex multivariable *-polynomial $V$. 

Indeed, they were able to establish that such an SDE has a unique stationary distribution $\mu_V$ satisfying the Schwinger-Dyson equation
$$ \mu_V \otimes \mu_V (\partial_i P) = \mu_V (D_i VP)$$
where $\partial_i$ is the non-commutative partial difference quotient and $D_i$ is the cyclic partial derivative.
By using the fact that they also had convergence in norm to this distribution from all initial data, they were able to show that the von Neumann algebra $M_V$ generated by operators with joint law $\mu_V$ is a factor with the Haagerup property. They also proved that $M_V$ has finite free entropy dimension and hence is prime and has no Cartan subalgebras. All of this provided evidence for the conjecture of Voiculescu that $M_V$ is isomorphic to a free group factor.

Recall that a von Neumann algebra $M \subseteq \mathbb{B}(H)$ is said to have property AO if there are ultraweakly dense C$^*$ subalgebras $A \subseteq M$ and $B \subseteq M'$ with $A$ locally reflexive and such that the  *-homomorphism $\Phi: A \otimes B \rightarrow \mathbb{B}(H)/\mathbb{K}$ given by
$$\Phi\left(\sum a_i \otimes b_i\right) = \pi\left(\sum a_i b_i \right)$$
is continuous with respect to the minimal tensor norm.

In this short paper we will demonstrate how the above techniques can be used to prove that $M_V$ has this property AO of Ozawa, and is thus solid by Theorem 6 in \cite{ozsolid}. This result adds to the evidence of the conjecture above as solidity (and AO) are well-known properties of free group factors.

\section{Preliminaries}

Denote the set of polynomials in $m$ non-commuting indeterminates $\left(X_1,\dots,X_m\right)$ by $\cxm$. Consider this as a subalgebra of $\cxmnsa$, the *-algebra of polynomials in the $m$ indeterminates and their formal adjoints. We will say that a polynomial $V \in \cxm$ is self-adjoint if $V(X_1,\dots,X_m)^* = V(X_1^*,\dots,X_m^*)$.

Recall the cyclic gradient $D$ of Rota, Sagan, and Stein, which is linear and given on any noncommutative multinomial $P \in \cxm$ by $DP = (D_1P,\dots,D_mP)$ with
$$D_i P = \sum_{P=Q X_i R} RQ.$$

Next, recall the non-commutative difference quotient $\partial: \cxm \rightarrow \cxm \otimes \cxm$ which is again linear and given on a multinomial $P$ by $\partial P = (\partial_1 P,\dots, \partial_m P)$ with
$$\partial_i P = \sum_{P=Q X_i R} Q \otimes R.$$

Finally, define for a pair of $m$-tuples of elements $X=(X_1,\dots,X_m)$ and $Y=(Y_1,\dots,Y_m)$ in any *-algebra the notation
$$X.Y = \frac{1}{2} \sum_{i=1}^m (X_iY_i^* +  Y_iX_i^*).$$

Then, for $V \in \cxm$, we say that $V$ is $(c,M)$-convex  if for any $m$-tuples $X=(X_1,\dots,X_m)$ and $Y=(Y_1,\dots,Y_m)$ of operators in any $C^*$ algebra $A$ satisfying $\Vert X \Vert_A , \Vert Y \Vert_A \le M$ we have
$$\left(DV(X) - DV(Y)\right).(X-Y) \ge c(X-Y).(X-Y).$$
(This is to be understood as an operator inequality in $A$.)

We will be considering the solutions to the multivariable SDE
$$ dX_t = dS_t - \frac{1}{2}DV(X_t)dt $$
for $V$ a $(c,M)$-convex polynomial. As such, we will fix an ambient free probablility space $(\mathcal{A},\tau)$ generated by a free Brownian motion $S$. For more information on SDE and free Brownian motion, see \cite{BS,BS2}.

Let $V$ be a (c,M)-convex polynomial in $\cxm$ and recall the following result (\cite{gushdiffusions} Lemma 2.1):

\begin{thm}\label{exist}
There exist finite constants
$$M_0=M_0(c, \|DV(0).DV(0)\|),$$
$$B_0=B_0(c, \|DV(0).DV(0)\|),$$
and 
$$b=b (c,\|DV(0).DV(0)\|, M) \ge B_0$$
so that whenever $M\geq M_0$  and $Z$ is an $m$-tuple with $\Vert Z \Vert < b$, there exists a unique solution  $X^Z_t$ to the SDE
$$ dX^Z_t = dS_t - \frac{1}{2}DV(X^Z_t)dt, \qquad t\in [0,+\infty) $$ with the initial data $X_0^Z=Z$.

Moreover, in this case,

\begin{eqnarray*}
\Vert X^Z_t \Vert \leq M,&\qquad& \forall\, t\in [0,+\infty),\\
\limsup_{t\to\infty} \Vert X^Z_t \Vert \leq B_0,&& \\
X^Z_t \in C^*(Z, S_q : q\in [0,t]),&\qquad&\forall\, t\in[0,+\infty). 
\end{eqnarray*} 

If $V$ is self-adjoint  $(c,M)$-convex and $X_0=Z$ is self-adjoint, then the above results hold and additionally, $X_t$ remains self-adjoint for all $t\ge 0$.
\end{thm}

Thus, if $V$ is sufficiently locally convex, and our initial data is appropriately bounded, we then have a unique bounded solution that exists for all time to our desired free diffusion equation. Now, let $V$ be as above and assume additionally that $V$ is self-adjoint. If we change the initial data for our SDE, we have the following asymptotic uniqueness result (\cite{gushdiffusions} Theorem 2.2):

\begin{thm}\label{converge}
Let $M_0$, $B_0$ and $b$ be as in Theorem~\ref{exist}, and assume that
$M\geq M_0$, and that $Z$ is an 
$m$-tuple of operators with $\Vert Z\Vert_\infty < b.$ 
Consider the unique solutions $X_t^Z$, $X_t^0$ to the free SDE
$$dX_t = dS_t - \frac{1}{2} DV(X_t)dt $$

with initial conditions $X_0^Z=Z$, and $X_0^0 = 0$ respectively.  
Then
\renewcommand{\theenumi}{\roman{enumi}}
\begin{enumerate}
\item $\Vert X_t^Z - X_t^0 \Vert_\infty \to 0\qquad\textrm{as $t\to\infty$}.$

\item The law of $(X_t^Z)$
converges to a stationary 
law  $\mu_V\in\cxm'$
which satisfies for all $P\in\cxm$ 

\begin{equation}\label{SD}
\sum_{i=1}^m\mu_V\otimes\mu_V(  \partial_i D_i P)=\sum_{i=1}^m \mu_V(D_i V D_i P).
\end{equation}

Moreover, for all $i\in\{1,\ldots,m\}$,
$$\mu_V(X_i^k) \le B_0^k.$$

Any law $\nu\in\cxm'$ of variables bounded in operator norm by $b$
that satisfies \eqref{SD} is such that $\nu= \mu_V$.
\end{enumerate}
\end{thm}

\section{Result}
Assume that $V$ is a self-adjoint $(c,M)$-convex polynomial in $m$ variables with $M \ge M_0$ for $M_0$ as above.

\begin{thm}\label{mainresult}
Let $Z$ be an $m$-tuple of self-adjoint operators satisfying $\Vert Z \Vert_\infty < b$ and  having the stationary law $\mu_V$ from above. Then $W^*(Z)$ has Ozawa's property AO, and is hence solid. 
\end{thm}

We first establish some notation and prove two lemmas.

Let $N=W^*(Z)$, $A=C^*(Z)$, and $M=W^*(Z,S_t: t \ge 0)$. Endow $M$ with the trace $\tau$ from the ambient free probability space, and let $H = L^2(M,\tau)$. Denote the unique solutions to the free SDE
$$dX_t = dS_t - \frac{1}{2} DV(X_t)dt $$
with initial data $0$ and $Z$ by $X^0_t$ and $X^Z_t$, respectively.
As the law of $X^Z_t$ is stationary, we have, for each $t$, homomorphic embeddings $\theta_t: N \rightarrow M$ satisfying $\theta_t(Z)=X^Z_t.$ Let $B=J_H A J_H$, and define $\tilde\theta_t(J_H a J_H) = J_H\theta_t(a)J_H$ for $a \in A$.

\begin{lem}\label{polyEst}
Let $\sum_{i=1}^n a_i \otimes b_i$ be an element in the algebraic tensor product $A \otimes B$. For every $\epsilon > 0$ there exists $p_i, q_i \in \cxmnsa$ and $t_0 \ge 0$ such that for all $t \ge t_0$,
$$ \left\Vert \sum_{i=1}^n \theta_t(a_i) \otimes \tilde\theta_t(b_i) - 
\sum_{i=1}^n p_i(X^0_t) \otimes J_Hq_i(X^0_t)J_H \right\Vert_{\text{min}} < \epsilon $$
and
$$ \left\Vert \sum_{i=1}^n \theta_t(a_i) \tilde\theta_t(b_i) - 
\sum_{i=1}^n p_i(X^0_t) J_Hq_i(X^0_t)J_H \right\Vert_\infty < \epsilon$$
where $\Vert \cdot \Vert_{\text{min}}$ denotes the minimal tensor norm on $M \otimes J_HMJ_H$.
\end{lem}

Proof:
First note that if $\delta > 0$ and $a \in A$, then we can find a $p \in \cxmnsa$ and a $t_0 \ge 0$ such that for all $t \ge t_0$,
$$\Vert \theta_t(a) - p(X^0_t)\Vert_\infty < \delta.$$
Indeed, first choose $p$ so that 
$$\Vert a - p(Z)\Vert_\infty < \delta/2.$$ 
Then, as by Theorem~\ref{converge} we have that $\Vert X^Z_t - X^0_t \Vert_\infty \rightarrow 0$ as $t \rightarrow \infty$, choose $t_0$ so that for all $t \ge t_0$,
$$\Vert \theta_t(p(Z)) - p(X^0_t)\Vert_\infty = \Vert p(X^Z_t) - p(X^0_t) \Vert_\infty < \delta/2.$$
The triangle inequality then implies that these are the desired $p$ and $t_0$. 

Note also that conjugation by $J_H$ shows that we can obtain a similar corresponding statement for approximating $\tilde\theta_t(b)$ by $J_H q(X^0_t) J_H$ for $q \in \cxmnsa$ and $t$ sufficiently large.

Now, fix $1 \ge \epsilon > 0$, and for each $i$ choose $p_i, q_i \in \cxmnsa$ and $t$ such that for all $t \ge t_0$,
$$\Vert \theta_t(a_i) - p_i(X^0_t) \Vert_\infty 
< \frac{\epsilon}{2n(\Vert b_i \Vert_\infty + 1)}$$
and
$$\Vert \tilde\theta_t(b_i) - J_H q_i(X^0_t) J_H \Vert_\infty
< \frac{\epsilon}{2n(\Vert a_i \Vert_\infty + 1)}.$$
Then we have by a simple estimate that
$$\Vert \theta_t(a_i) \otimes \tilde\theta_t(b_i) - p_i(X^0_t) \otimes J_Hq_i(X^0_t)J_H \Vert_\text{min}
< \epsilon/n,$$
and
$$\Vert \theta_t(a_i) \tilde\theta_t(b_i) - p_i(X^0_t) J_Hq_i(X^0_t)J_H \Vert_\infty
< \epsilon/n.$$
Thus the triangle inequality again implies the lemma. $\square$

\begin{lem}\label{cptEst}
There exists $\alpha > 0$ such that for any $p_i,q_i \in \cxmnsa$ where $i\in\{1,\dots,n\}$, any $\epsilon > 0$, and any $t \ge 0$ there exists a compact operator $T$ such that
$$ \left\Vert \sum_{i=1}^n p_i(X^0_t)J_H q_i(X^0_t) J_H + T \right\Vert_\infty \le
\alpha \left\Vert \sum_{i=1}^n p_i(X^0_t) \otimes J_H q_i(X^0_t) J_H \right\Vert_\text{min} + \epsilon.
$$
\end{lem}

Proof: Let $C = C^*(S_t: t \ge 0)$, $D = J_H C^*(S_t: t \ge 0) J_H$, and note that $C$ is nuclear. Thus the *-homomorphism $\Phi: C \otimes D \rightarrow \mathbb{B}(H)/\mathbb{K}$ given by
$$\Phi\left(\sum c_i \otimes d_i\right) = \pi\left(\sum c_i d_i \right)$$
(for $\pi$ the canonical homomorphism into the Calkin algebra)  is continuous with respect to the minimal tensor norm on $C \otimes D$.

So, there exists an $\alpha > 0$ such that for every $c_1, \dots c_n \in C$, $d_1, \dots, d_n \in D$ and $\epsilon > 0$ there exists a compact $T_\epsilon$ such that
$$\left\Vert \sum_{i=1}^n c_i d_i + T_\epsilon \right\Vert_\infty \le 
\alpha \left\Vert \sum_{i=1}^n c_i \otimes d_i \right\Vert_{\text{min}} +\epsilon$$

Note that by Theorem~\ref{exist} we have that $X^0_t \in C$ for any $t$, and so letting $c_i=p_i(X^0_t)$, $d_i = J_Hq_i(X^0_t)J_H$ and setting $T=T_\epsilon$ proves the lemma. $\square$
\vspace{12 pt}

Proof of Theorem \ref{mainresult}: We will show that the *-homomorphism $\Psi: A \otimes B \rightarrow \mathbb{B}(L^2(N,\tau))/\mathbb{K}$ given by
$$\Psi\left(\sum a_i \otimes b_i\right) = \pi\left(\sum a_i b_i \right)$$
is continuous with respect to the minimal tensor norm on $A \otimes B$. (Here we have identified $B$ with its restriction to $L^2(N,\tau) \subset H$.) Note that this indeed suffices as by Lemma 4.3 in \cite{gushdiffusions}, $A$ is exact and hence locally reflexive.

To this end, let $\epsilon>0$ and let $\alpha$ be as in Lemma~\ref{cptEst}. By Lemma~\ref{polyEst} choose $p_i$, $q_i$ and $t$ such that
$$ \left\Vert \sum_{i=1}^n \theta_t(a_i) \otimes \tilde\theta_t(b_i) - 
\sum_{i=1}^n p_i(X^0_t) \otimes J_Hq_i(X^0_t)J_H \right\Vert_{\text{min}} < \frac{\epsilon}{2(\alpha+1)} $$
and
$$ \left\Vert \sum_{i=1}^n \theta_t(a_i) \tilde\theta_t(b_i) - 
\sum_{i=1}^n p_i(X^0_t) J_Hq_i(X^0_t)J_H \right\Vert_\infty < \frac{\epsilon}{2(\alpha+1)}. $$

Then, apply Lemma~\ref{cptEst} to find a compact operator $T$ such that

$$ \left\Vert \sum_{i=1}^n p_i(X^0_t) J_Hq_i(X^0_t)J_H + T \right\Vert_\infty \le
\alpha \left\Vert \sum_{i=1}^n p_i(X^0_t) \otimes J_Hq_i(X^0_t)J_H \right\Vert_{\text{min}} + \epsilon/2.$$

We then have by another simple estimate that
$$\left\Vert \sum_{i=1}^n \theta_t(a_i) \tilde\theta_t(b_i) + T \right\Vert_\infty \le
\alpha \left\Vert  \sum_{i=1}^n \theta_t(a_i) \otimes \tilde\theta_t(b_i) \right\Vert_{\text{min}} + \epsilon.$$
If we restrict the operator on the left-hand side to $L^2(\theta_t(N),\tau)$, we will obtain the inequality
$$\left\Vert \sum_{i=1}^n \theta_t(a_i) \tilde\theta_t(b_i) + S \right\Vert_\infty \le
\alpha \left\Vert  \sum_{i=1}^n \theta_t(a_i) \otimes \tilde\theta_t(b_i) \right\Vert_{\text{min}} + \epsilon$$
for $S = e_N T e_N$ a compact operator on $L^2(\theta_t(N),\tau)$. By identifying N with $\theta_t(N)$ and $L^2(N,\tau)$ with $L^2(\theta_t(N),\tau)$, we thus get a compact operator $R$ on $L^2(N,\tau)$
$$\left\Vert \sum_{i=1}^n a_i b_i + R \right\Vert_\infty \le
\alpha \left\Vert  \sum_{i=1}^n a_i \otimes b_i \right\Vert_{\text{min}} + \epsilon$$
and this proves the theorem. $\square$

\end{document}